\documentclass[12pt,a4paper]{article}

\usepackage{amsmath}
\usepackage{amssymb}
\usepackage{latexsym}
\usepackage{graphics}
\usepackage{graphicx}
\usepackage[ruled,vlined]{algorithm2e}

\newcommand{\co}{{\mathbb C}}

\newcommand{\re}{{\mathbb R}}
\newcommand{\n}{{\mathbb N}}

\newcommand{\cA}{{\cal{A}}}

\newcommand{\cB}{{\cal{B}}}
\newcommand{\cD}{{\cal{D}}}
\newcommand{\cI}{{\cal{I}}}

\newcommand{\RP}{{\cal{RP}}}

\newcommand{\cP}{\cal P}

\newcommand{\be}{{\boldsymbol{e}}}
\newcommand{\bh}{{\boldsymbol{h}}}
\newcommand{\by}{{\boldsymbol{y}}}
\newcommand{\bx}{{\boldsymbol{x}}}
\newcommand{\bz}{{\boldsymbol{z}}}

\newtheorem{theorem}{Theorem}
\newtheorem{prop}{Proposition}
\newtheorem{lemma}{Lemma}
\newtheorem{cor}{Corollary}
\newtheorem{remark}{Remark}

\newtheorem{defi}{Definition}
\newtheorem{conj}{Conjecture}

\date{}

\begin{document}

\author{Vladimir Yu. Protasov
\thanks{University of L'Aquila (Italy),  {e-mail: \tt\small
vladimir.protasov@univaq.it}}}

\title{Generalized Markov-Bernstein inequalities\\ 
and stability of dynamical systems
\thanks{
The research is supported by the RFBR grant 20-01-00469.
}}

\maketitle

\begin{abstract}

The Markov-Bernstein type inequalities between the norms of functions and of their derivatives 
are analysed for complex exponential polynomials. We establish a relation between the sharp constants in those inequalities and the stability problem for linear switching systems.  
In particular, the maximal discretization step is estimated. We prove the monotonicity of the sharp constants with respect to the exponents, provided those exponents are real. 
This gives asymptotically tight uniform bounds and the general form of the extremal  polynomial. 
The case of complex exponent is left as an open problem.

\bigskip

\smallskip

\begin{flushright}
\noindent  \textbf{AMS} {\em    41A50,    93D20, 41A17, 37N35}

\end{flushright}

\end{abstract}
\bigskip

\begin{center}
\textbf{\large{1. Introduction}}
\end{center}
\medskip

The Markov-Bernstein type inequalities associate the norm of a polynomial 
to the norm of its derivative of $\ell$-th order. The Bernstein  inequality 
states that for a trigonometric polynomial of degree~$n$, it holds that~$\|p'\| \, \le \, n\, \|p\|$, where  $\|p\| = \sup_{t \in \mathbb{T}}|p(t)|$ is 
the uniform norm on the period. The Markov inequality states that 
for an algebraic polynomial of degree~$n$ on the segment~$[-1,1]$, it holds~$\|p'\| \, \le \, n^2\, \|p\|$, and the equality is attained 
precisely when 
$p$ is proportional to the Chebyshev polynomial~$T_n$. 
We are interested in generalizations of those inequalities 
to the exponential polynomials, i.e., to the functions of the form
 $p(t)\, = \, \sum_{k=1}^n p_k e^{-h_k t}$, where 
 $h_1, \ldots , h_n$ are given complex numbers 
(in case of multiplicity~$m$ the corresponding exponent is multiplied by the 
powers~$t^k, \, k = 0, \ldots , m-1$). 
For every~$\ell \ge 1$, for 
arbitrary~$h_k$ and for a segment  (domain), we consider the  
inequality 
$\|p^{(\ell)}\| \, \le \, C\, \|p\|$ 
with the sharp constant~$C$.  The case when~$\ell = 1$,  
the exponents are purely imaginary and form an arithmetic progression, 
corresponds (after a proper change of variables) to the Bernstein inequality.   
If the exponents are real and form an arithmetic progression, then this case corresponds to  the Markov inequality.

Note that in the case of real rational powers~$h_k$ 
the exponential polynomial~$p(t)$  
is transferred by the change of variables  $x=e^{-t/N}$  to  
an algebraic polynomial~$P(x)$. However, this does not solve the problem 
of finding the sharp constant, since the polynomial~$P(x)$
is not arbitrary: its degree can significantly exceed the number of its 
nonzero coefficients. Therefore, the constant 
obtained this way can be much  bigger than the sharp constant. 

We consider inequalities for derivatives on the positive half-axis~$\re_+$ 
naturally assuming that  ${\rm Im}\, h_k \, > \, 0$ and, therefore, 
all the polynomials~$p(t)$ tend to zero as $t \to +\infty$. 
The main problem is to find the precise value or at least to estimate 
the sharp constant~$C$ in the inequality 
$\|p^{(\ell)}\| \le C\, \|p\|$. For the  derivative~$p^{(\ell)}$, we always choose 
the norm in~$C(\re_+)$, while the norm~$\|p\|$ of the polynomial 
can be taken more general. In particular, a part of the results 
are obtained for an arbitrary monotone norm on~$\re_+$ 
(the monotonicity means that  $\|f\| \ge \|g\|$ if $|f(t)| \ge |g(t)|$ 
at every point~$\ t \in \re_+$). 

We are interested in the values of the constants~$C = C(\ell, \bh)$ 
in the Markov-Bernstein inequality for fixed exponents~$\bh$ 
and in the uniform constants for all values of~$\bh$
from a given domain. This problem is  motivated 
by applications in the study of trajectories of linear systems. 

\smallskip   

\textbf{{\em The problem of the maximal initial velocity of a bounded trajectory.}}
If a function 
$x \in C^n[0, +\infty)$ is a solution of a linear ODE with constant coefficients 
$x^{(n)} \, = \, \sum_{k=0}^{n-1} a_k x^{(k)}(t)$, then how large its initial 
 derivative~$x'(0)$ can be under the condition $\|x\|_{C[0, +\infty)} \le R$, 
 where 
$R> 0$ is given? In other words, the trajectory never leaves the ball 
of radius~$R$.  The same question can be formulated for the 
higher order derivatives~$x^{(\ell)}(0)$. 
The answer is that $|x^{(\ell)}(0)| \, \le \, R\, C(\ell, \bh)$, 
where  
$C(\ell, \bh)$ is the sharp constant in the Markov-Bernstein inequality for 
exponential polynomials over 
the system~$\{e^{-h_k t}\}_{k=1}^n$, where $\{-h_k \}_{k=1}^n$ 
are the roots of the characteristic polynomial of the ODE. 
Similar problems are considered for other norms, for example, for~$\|x\|_{L_p} \le R$. 

\smallskip   

\textbf{{\em The stability problem for linear switching systems}}. 
This is probably the most important application and the main motivation of this research. 
We realize it in Section~7. The linear switching system is a linear ODE  $\dot \bx (t) \, = \,  A(t)\bx(t), \ \bx(0) = \bx_0$, 
where~$\bx(t) \in \re^n$ and   $A(\cdot )$ is an arbitrary measurable function 
taking values on some compact set of $n\times n$ matrices~$\cA$.  
Such systems regularly arise in engineering applications. Their systematic study 
began in late 70-s, see~\cite{O, MP1, MP2, L} and references therein. 
A system is asymptotically stable if its trajectories $\bx(t)$ 
tend to zero as  $t \to +\infty$ for all functions~$A(\cdot )$. 
One of the methods of proving stability is the discretization of the system, i.e., 
its approximation by the discrete-time system: 
$\bx(k+1) \, = \, \bigl(I \, + \, \tau A(k)\bigr)\bx(k), \ k = 0, 1, \ldots $, 
where $I$ is the identity matrix. It is known~\cite{MP2}  that
if the system becomes stable after the discretization with 
some step~$\tau$, then it does  
for all smaller steps and it is stable (as a continuous-time system). 
The problem consists of finding the longest possible step~$\tau$ for which 
the converse is also true with some precision~$\varepsilon > 0$. This means that 
if the system $\dot \bx  \, = \,  A\bx$
is stable, then the discrete system  $\bx(k+1) \, = \, \bigl(I \, + \, 
\tau \, ( A(k) \, - \, 	\varepsilon I )\bigr)\, \bx(k)$ is stable as well. 
Since there are efficient methods for deciding the discrete-time stability~\cite{G, GP13, PJB, Mej, MR14}, it follows that an efficient estimate of the discretization step~$\tau$ 
makes those methods applicable also for continuous-time systems. 
This idea was developed in~\cite{BCM, BS, PJ1, PJ2}, etc.  
In Section~7 we estimate~$\tau$ in terms of the sharp constant~$C(\ell, \bh)$ 
in the Markov-Bernstein inequality for~$\ell = 2$ and $\bh = - {\rm sp}\,(A), \ A \in \cA$,
where ${\rm sp}\,(A)$ denotes the set of eigenvalues of the matrix~$A$. 
Those results are formulated in Theorems~\ref{th3} and~\ref{th5}. 
\smallskip   

\textbf{{\em The fundamental results}}. We obtain the estimates for the step size~$\tau$ 
by the constants~$C(\ell, \bh)$ for arbitrary vectors~$\bh$, while 
the estimates and the sharp values of these constants are found only for the 
real vectors. In the latter case the system of exponents is a Chebyshev system 
and, therefore, for finding the extremal function we can apply the 
alternance. A similar idea for estimating the constants in the Markov-Bernstein inequality 
for real exponential polynomials have been exploited in~\cite{BE0} -- \cite{BE3}, \cite{MN, New, R1}. We establish the comparison theorem  (Theorem~\ref{th10})
according to which the constant~$C(\ell, \bh)$ strictly decreases in
each component~$h_k$. This is true not only for a uniform norm in~$C(\re_+)$
but also for each $L_p$-norm, and, more generally, for all monotone norms~$\re_+$. 
Therefore, this constant reaches its  biggest value on the set $h_k \le \alpha_k, \ k = 1, \ldots , n$, at the point~$\bh = (\alpha_1, \ldots , \alpha_n)$. 
This result allows us to obtain uniform estimates for the constant~$C(\ell, \bh)$ 
for various ranges  of the vector~$\bh$. 
For example, under the conditions $h_k \le 1, \  k = 1, \ldots , n$, 
the maximal values are attained for the polynomial~$p(t) = e^{-t}R_{n-1}(t)$, 
where $R_{n-1}$ is the Chebyshev polynomial with the Laguerre 
weight~$e^{-t}$ (Theorem~\ref{th25}). Applying known estimates on the derivative of those polynomials~\cite{CLM, Fre,  MN, Sz} and especially the asymptotically sharp bound 
of V.Sklyarov~\cite{S}, we obtain uniform estimates 
for the constant in the Markov-Bernstein inequality (Theorem~\ref{th35}). 
They lead to practically applicable estimates for the discretization step~$\tau$ 
in the problem of stability of linear switching systems. 
Those estimates, however, can be used only for systems of real exponents, i.e., 
for switching systems defined by matrices with purely real spectra. 
\smallskip   

\textbf{{\em The problem of generalization to arbitrary exponents}}. 
A significant disadvantage of the obtained results is that they hold 
for real exponents. This restriction looks especially strange 
in the problem of the discretization steplength, where our estimates 
are true only  when  all the matrices of the system have real spectrum. 
This phenomenon does not cause problems for a concrete switching system since 
the value $C(\ell, \bh)$ can be found for every~$\bh$ as
the solution of a convex problem $p^{(\ell )}(0) \to \max, \ \|p\| \le 1$. 
However, to derive uniform estimates one needs the comparison theorem, whose proof  
cannot be extended to arbitrary systems of exponents since they do not form 
a Chebyshev system. We leave the proof of the comparison theorem 
for complex exponents as an open problem (Conjecture~1). If the answer is affirmative, 
then all estimates obtained for real exponents are extended to the complex case. 
\smallskip   

\textbf{{\em Notation}}.
We denote vectors by bold letters and scalars by standard letters. Thus, 
$\bh = (h_1, \ldots , h_n)$.  
We use the standard notation 
$\re_{+}$ and $\re_{++}$ for the set of nonnegative numbers and of positive numbers 
respectively; similarly  $\co_{++} = \bigl\{ z \in \co \ | \ {\rm Re}\,z \,  > \, 0\bigr\}$ 
is an open right half-plane of the complex plane;  
$L_p$-norm on~$\re_+$ will be denoted as $\|\cdot \|_p$.
In particular, if   $f \in C(\re_+)$, then $\|f \|_{\infty} = \|f\|_{C(\re_+)}$.

\bigskip

\begin{center}
{\large \textbf{2. Statement of the problem}}
\end{center}
\bigskip

For a given vector   $\bh = (h_1, \ldots , h_n) \in \co^n_{++}$, we 
consider the system of exponent 
 $E_{\bh} = \{e^{-h_k t}\}_{k=1}^n$. Some of the numbers~$h_k$ can coincide, 
 in this case the corresponding exponents are multiplied by powers of~$t$. 
 If, for instance, the components   $h_1, \ldots , 
 h_{r}$ are equal and are different from all others, i.e., 
 the exponent~$h_1$ has multiplicity~$r$, 
then the functions $e^{-h_{1}t}, 
 \ldots , e^{-h_{r} t}$
are replaced by $e^{-h_{1}t}, te^{-h_{1} t}, \ldots , t^{r-1}e^{-h_{1} t}$ respectively.
 A {\em polynomial} by the system~$\{e^{-h_k t}\}_{k=1}^n$, or a  {\em quasipolynomial}
is an arbitrary linear combination of those exponents with complex coefficients. 
The linear space by a given system on the half-line  $\re_+$ 
is denoted by  $\cP_{\bh}$. 
 This is an   $n$-dimensional subspace  
$C_0(\re_+)$ of functions that are continuous on~$\re_+$ 
an tend to zero as~$t \to +\infty$.  
The map $\bh \mapsto \cP_{\bh}$ is well-defined and continuous~\cite{K, KN}.
\smallskip

The real part of a quasipolynomial is a real linear combination of the 
functions $t^m e^{- \alpha_k t}\cos \beta_k t$, $\ 
t^m e^{- \alpha_k t}\sin \beta_k t, \ k = 1, \ldots , n$, where 
$\alpha_k, \beta_k$ are the real and the imaginary part of~$h_k$, and the degree~$m$
does not exceed the multiplicity of~$h_k$.  Linear combinations of those functions 
with real coefficients form the  {\em space of real quasipolynomials~$\RP_{\bh}$}.

Consider an arbitrary  {\em monotone norm}  $\|\cdot \|$ in the space~$\cP_{\bh}$. 
The monotonicity means that if 
 $|f_1(t)| \ge |f_2(t)|$ for all~$t \ge 0$, then $\|f_1\| \ge \|f_2\|$.  
 For example, the $L_p$-norms, $p \in [1, +\infty]$, the weighted 
 $L_p$-norms are monotone. For every monotone norm and for an arbitrary functional 
  $F: {\cP_{\bh}} \to \re$, we consider the following problem:
  \begin{equation}\label{markov2}
\left\{
\begin{array}{l}
F(p) \ \to \ \min \\
p \in {\cP_{\bh}}\, , \ \|  p \| \le 1\, .
\end{array}
\right.
\end{equation}
The value of this problem (the minimal value of the objective function)
will be denoted by~$\Phi(\bh)$.
Since this problem is convex,~$\Phi(\bh)$ can be found by standard tools of 
convex programming. We, however, are interested not in the 
numerical solution  but in the description of the extremal polynomial 
and in uniform estimates for~$\Phi(\bh)$ 
 over all vectors~$\bh$ from a given set. We deal with the sets~$\cD_n \ = \, \bigl\{\bz \in \co^n \ | \ |z_k| \le 1, \ 
{\rm Re}\, z_k \, > \, 0, \ 
\ k = 1, \ldots, n\bigr\}\, $ and $\, \cI_n \, = \, \bigl\{\bx \in \re^n \ | \ 
0< \, x_k \, \le 1, \ k = 1, \ldots, n \bigr\}$.  
We also use  simplified  notation: 
$\cD_1 = \cD$ (half-disc) and $\cI_1 = \cI$ (half-interval~$(0,1]$). 

\medskip

In this section we consider the problem~(\ref{markov2}) 
in the uniform norm~$\|\cdot \|_{\infty}$
on $\re_{+}$  
for the functional $F(p)\, = \, - \, p^{(\ell)}(0)$, where   $p^{(\ell)}$ 
is the $\ell$-th derivative of $p$, $\ell \ge 1$. 
Thus, one finds the value of the  $\ell$-th derivative of the polynomial~$p$ 
at zero under the assumption  $\|p\|_{\infty} \le 1$. 
This is equivalent to finding the maximal norm $\|p^{(\ell)}\|_{\infty}$ 
on the unit ball, which is equal to the sharp constant in the Markov-Bernstein inequality 
for polynomials from~$\cP_{\bh}$. 
\begin{defi}\label{d.15}
For a given vector~$\bh \in \co_{++}^n$, we set 
$$
M_{\, \ell}(\bh) \ = \ \max_{\|p\|_{\infty} \le 1, \ p  \in \cP_{\bh}}
\ \bigl\| \, p^{\, (\ell)}\, \bigr\|_{\infty}, 
$$
and 
$$
M_{\, \ell, n} \ = \ \max_{\bh \in \cD_{n} } M_{\ell}(\bh) \, .
$$
If we search the extremal polynomial among the real exponential 
polynomials with real coefficients~$h_i$, then 
$\bh$ is the real positive vector and the set $\cD_n$ 
is the non-closed unit cube
$\cI_n \, = \, (0,1]^n$:  
$$
m_{\, \ell}(\bh) \ = \ \max_{\|p\|_{\infty} \le 1, \ p  \in \RP_{\bh}}
\ \bigl\| \, p^{\, (\ell)}\, \bigr\|_{\infty}, \qquad 
m_{\, \ell, n} \ = \ \max_{\bh \in \cI_{n} }  m_{\, \ell}(\bh) \, .
$$
\end{defi}
Clearly, 
$m_{\, \ell, n} \, \le \, M_{\, \ell, n}$. 
Our main conjecture is that those two numbers are actually equal: 
\begin{conj}\label{conj10}
For all $\ell , n \in \n$, we have  $m_{\, \ell, n} \, = \, M_{\, \ell, n}$. 
\end{conj} 
The value   $M_{\ell}(\bh)$ is the sharp constant in the Markov-Bernstein type inequality 
for exponential polynomials: 
$\|p^{(\ell)}\|_{\infty} \ \le \ M_{\ell}(\bh)\, \|p\|_{\infty}, \ p\in \cP_{\bh}$. 
In fact, the same constant  $M_{\ell}(\bh)$ is sharp in the same inequality 
also for the class of real exponential polynomials~$\RP_{\bh}$. 
For the proof, it suffices to observe that by multiplying 
the polynomial~$p$ by a proper number~$\, e^{i\alpha}, \ \alpha \in \re$,   
and by translating we obtain the equality~$\|p^{(\ell)}\|_{\infty} \, = \, {\rm Re}\, p^{(\ell)}(0)$. 
Thus, $M_{\ell}( \bh)$ is the value of the problem~${\rm Re} \, p^{(\ell)}(0) \to \max\, , \ 
\|p\|_{\infty}\le 1$.  The replacement of the polynomial $p$ by its real part, which is a real 
quasipolynomial~$q$, does not change the value ${\rm Re} \, p^{(\ell)}(0) \, = \, q^{(\ell)}(0)$ and does not increase the norm. Therefore, the restriction of our problem to the class of real 
quasipolynomials does not change its value. 
\smallskip 

Then we will return to an arbitrary monotone norm in the problem~(\ref{markov2})
and will focus on real exponents~$h_1, \ldots , h_n$. 
We will prove the comparison theorem  (theorem~\ref{th10})
according to which the value of the problem strictly increases in each component~$h_i$. 
This yields that the  value~$m_{\ell}( \bh)$ also increases in each component~$h_i$. 
 This allows us to find the optimal value 
 over various domains of the parameter~$\bh$ and obtain uniform bounds over those 
 domains. Further we apply those estimates to the 
 problem of stability
  of linear switching systems.

\bigskip

\begin{center}
\large{ \textbf{3. The case of real exponents}}
\end{center}
\vspace{1cm}

A system  $E_{\bh} = \{e^{-h_k t}\}_{k=1}^n$ is called real if  
 $h_k$ are all real. In this case all those numbers are strictly positive and 
 are assumed to be arranged in ascending order: $0< h_1 \le \cdots \le h_n$. 
The space of polynomials by the system~$E_{\bh}$ on $\re_+$ with real coefficients 
is called the space of real exponential polynomials  (or real quasipolynomials) 
and is denoted by~$\RP(\bh)$. The main results of this section are easily generalized 
to polynomials with complex coefficients (but still with real exponents~$h_i$). 
For the sake of simplicity, we consider only the case of real coefficients. 
We need several basic properties of the space~$\RP(\bh)$.
\bigskip

1) $E_{\bh}$ is a Cartesian system, i.e., 
it satisfies the Cartesian rule of counting zeros, see, for example,~\cite{Dz}. 
 This implies in particular  that if the polynomial  $p \in \cP_{\bh}$ has $n-1$ 
zeros on~$\re_+$ counting multiplicities, then all its  
coefficients are nonzero and their signs alternate.
\smallskip
 
 The Cartesian property can be proved by approximating all the numbers 
 $h_i$ by rational numbers   $ \tilde h_i$ and by the change of variables 
  $x = e^{-t/N}$, where   $N$ is such that 
 $\tilde h_i N \in \n\, , \, i = 1, \ldots , n$.
Thus, $E_{\bh}$ is approximated with arbitrary precision 
by algebraic polynomials, which, as we know, form a  Cartesian system. 
\medskip

 2) Every Cartesian system is also a Chebyshev system 
(or Haar system), i.e., every nontrivial polynomial from~$\cP_{\bh}$
has at most   $n-1$ zeros~\cite{KS, KN, Dz}. Thus, $E_{\bh}$ 
is a Chebyshev system on~$\re_+$. 
For every~$n$, for arbitrary~$t_i \in \re^+$, and for arbitrary numbers 
$c_i \in \re, \, i = 1, \ldots , n$,  there exists a polynomial~$p \in \cP_{\bh}$
for which~$p(t_i) = c_i$.

 \bigskip

 3) For each~$\ell \in \n$, the $\ell$-th derivative of the polynomial  $p \in \cP_{\bh}$ 
has at most  $n-\ell$  zeros on 
 $\re_+$ counting multiplicities.
\bigskip

Indeed, if we denote all zeros of $p$ by $\, a_1 < \ldots < a_k < +\infty$, 
assuming now that they are all different, and adding an extra zero
$a_{k+1} = +\infty$  (since $p(+\infty) = 0$)
we come to the conclusion that each interval 
 $(a_i, a_{i+1}), \, \, i = 1, \ldots , k$, contains at least one zero of the 
 derivative~$p\,'$.
Therefore, $p\,'$ has at least  $k$ 
zeros on $\re_+$. Applying induction we extend this property for all derivatives  $p^{(\ell)}$. 
The case of multiple roots follows by a limit passage.   
\bigskip

\begin{center}
\large{ \textbf{4. The comparison theorem}}
\end{center}
\bigskip

For a system of real exponents~$\bh$, it is possible 
not only to compute the value~$m_{\ell}(\bh)$  but also to 
analyse its behaviour as a function of the arguments~$h_1, \ldots , h_n$. 
In other words, one can find the asymptotics of the  $\ell$-th derivative  
of the polynomial in the unit ball of the space~$\cP_{\bh}$ depending on the vector~$\bh$.  
Moreover, this can be done not only for the unit ball in the uniform norm~$\|\cdot\|_{\infty}$, but also for every monotone norm in~$\re_+$, in particular, for the~$L_p$-norm. 
As we mentioned in the previous section, it is sufficient to solve this problem 
the space of quasipolynomials with real coefficients~$\RP_{\bh}$. 

Consider the problem~(\ref{markov2}) for the functional  $F(p) = - p^{(\ell)}(0)$
and for a fixed monotone norm~$\|\cdot \|$ in $\RP_{\bh}$. 
Thus, it suffices to find the maximum of the  $\ell$-th derivative 
at zero under the assumption~$\|\cdot \| \le 1$. Denote by $\Phi_{\ell}(\bh)$
the value of this problem. The polynomial~$p$
on which this maximum is attained has 
$n-1$ zeros (counting multiplicities) on the ray~$\re_+$. Otherwise, 
there would exist a polynomial 
 $q \in \cP_{\bh}$ that has the same roots as~$p$ does and such that 
  $p(t)q(t) < 0$ for all  $t \in \re_+$ different from the roots of~$p$
  and for which~$p^{(\ell)}(0)q^{(\ell)}(0) > 0$. 
Then for sufficiently small~$\lambda$, we have 
 $\|p+q\| < 1$  and  $F(p+q) > F(p)$, which contradicts to the optimality of~$p$.
\begin{theorem}\label{th10} (\textbf{comparison theorem})
Let $\|\cdot \|$ be an arbitrary monotone norm in~$C_0(\re_+)$ and 
$\Phi_{\ell}(\bh)$ be the maximal value of~$p^{(\ell)}(0)$ for all possible $p\in \RP_{\bh}$
such that~$\|p\|\le 1$. 
Then, if the vectors $\,  \bh', \bh\, \in \, \re_{++}$  are such that         
  $\,  \bh' \ge \bh\, $ and  $\, \bh' \ne \bh$, then for all $\ell \in \n$, 
  we have   $\ \Phi_{\ell}(\bh') > \Phi_{\ell}(\bh)$.
\end{theorem}
We see that the value of the problem~(\ref{markov2}) 
increases in every variable~$h_i$. 
Applying this theorem for the~$L_{\infty}$-norm, we obtain 
\begin{cor}\label{c15} 
If $\,  \bh' \ge \bh\, $ and $\, \bh' \ne \bh$, then for each $\ell \in \n$, we have   $\ m_{\ell}(\bh') \, > \, m_{\ell}(\bh)$.
\end{cor}
For~$\ell\, \, =1$, this corollary is analogous to the comparison theorem 
for hyperbolic sines~\cite{BE2}, but the method of the proof is 
different and is based on the following key fact: 
a small perturbation of a Chebyshev system   $E_{\bh}$ 
generates a bigger Chebyshev system (i.e., 
a system of a larger number of functions). This provides an additional 
degree of freedom for the choice of the parameter and makes it possible  to reduce the 
norm of the polynomial and to simultaneously increase the value of the objective function $F$.

We will use the following fact that is simply derived from the convexity of a norm. 
\begin{lemma}\label{l10}
If  $p$ and $q$ are elements of a normed space, $\|p\| = 1$ and 
$\|p +  q\| < 1$, then there exists a positive constant  $c$ such that  $\|p + \lambda q\|\, < \, 1\, - \, c \lambda$
for all $\, \lambda \in (0, 1]$.
\end{lemma}

{\tt Proof of Theorem~\ref{th10}.} Without loss of generality we assume that 
${p^{(\ell)}(0) > 0}$. It suffices to consider the case
when all coordinates of the vectors  $\bh$ and $\bh'$, except for one of them,
say, the $r$th one, coincide:  $\bh' = \bh + \tau \be_r\, , \ \tau > 0$. 
After this the theorem is proved by step-by-step changing the coordinates. 
This, in turn, is sufficient to prove for small variations of
$\tau$, because locally increasing functions increase globally. 
Further, it can be assumed that the exponent~$h_{\, r}$ is simple. 
Indeed, if it is multiple, then the inequality  $F(p_{\, \bh'}) \ge F(p_{\, \bh})$ 
follows from the continuity,
 and the strict inequality $F(p_{\, \bh'}) > F(p_{\, \bh})$ 
 follows by replacing the shift   $h_{\, r} \to h_{\, r}'$
by two consecutive shifts  $h_{\, r}\to h_{\, r}'' \to h_{\, r}'$ with simple 
exponent~$h_{\, r}''$. The first shift does not reduce $F$ and the second one 
increases it. Thus, the exponent $h_{\, r}$ is simple. Assume that so are all other $h_k$
 (in the general case the proof is similar).

The maximum of   $F$ is attained at some polynomial  $p(t)  = \sum_{k=1}^n p_k e^{-h_kt}$, 
that has 
 $n-1$ roots on~$\re_+$. Suppose all of them are simple: $0< \mu_1 < \ldots < \mu_{n-1} < \infty$. The general case is treated in the same way. 
Consider the vector $\by \in \re^n$ and a number $\delta > 0$
which will be specified later.
For arbitrary  $\lambda > 0$, take  a $(\lambda \by, \lambda \delta)$-variation 
of the polynomial~$p$:
$$
p_{\lambda}(t)\  = \
(p_r \, + \, \lambda y_r)\, e^{-(h_{\, r}+\lambda\delta)t}\, + \, \sum_{k\ne r} (p_k \, + \, \lambda y_k)\, e^{-h_kt}
$$
Thus, $p_{\lambda} \in \cP_{\bh'}$ with $\bh' = \bh + \lambda \delta \be_r$. 
We aim to find  $\by$ and $\delta$ such that  $p_{\lambda}$ has the small norm and 
a large value $F$ for sufficiently small $\lambda > 0$. 
Then  $e^{-(h_{ r}+\lambda\delta)t} \, = \, e^{-h_{ r} t} (1 - \lambda \delta t) \, + \, O(\lambda^2)$, and therefore 
\begin{equation}\label{ineq1}
p_{\lambda}(t) \, - \, p(t)\ = \  \lambda \, \Bigl[ \, - \, y_r \, \delta \, t e^{-h_{ r}t}\ + \ \sum_{k=1}^n y_{k} e^{-h_kt}\, \Bigr] \ + \
O(\lambda^2)\, , \qquad \lambda \to 0\, .
\end{equation}
Since   $\, te^{-h_{ r}t} \to 0$ as  $t \to +\infty$,
the latter term in~(\ref{ineq1}) is  $O(\lambda^2)$ uniformly for all $t \in [0, +\infty)$.
Hence, the polynomial  $p_{\lambda}- p$ can be 
approximated with precision  $O(\lambda^2)$ by a polynomial 
that consists of $n+1$ exponents (the left-hand side of~(\ref{ineq1})).
Now we choose the coefficients of those polynomials so that $\|p_{\lambda}\| < 1$ and  $F(p_{\lambda}) > F(p)$.
Consider the polynomial $q(t) = \sum_{k=1}^n q_ke^{-h_k t} \, + \, q_{r, 1}te^{-h_r t}$, 
that consists of   $n+1$ exponents ($h_{\, r}$ has multiplicity~$2$), 
and that vanishes in $n$ points  $-\alpha = \mu_0, \mu_1, \ldots , \mu_{n-1}$, where $\alpha < 0$, at all points  $t \in \re_+$ different from the roots~$\mu_k$
and, therefore,   $|p(t) + q(t)| < 1$, hence   $\|p+q\|< 1$. 
Furthermore, the derivative  $q^{(\ell)}(t)$ 
changes its sign  $n$ times on the interval  $[-\alpha, +\infty)$, 
while  $p^{(\ell)}(0)$ does 
  $n-1$ times and the signs of those polynomials are different as $t \to +\infty$. 
  Therefore, they have the same sign at the point  $t = -\alpha$. 
Thus, $q^{(\ell)}(-\alpha) > 0$. Choosing a small $\alpha$ we can assume that  $q^{(\ell)}(0) > 0$, and so $F(q) > 0$.

Now we choose the coefficients  $y_1, \ldots , y_n, \delta$ so that the polynomials 
in the right-hand side of~(\ref{ineq1}) will coincide with $q$. 
In this case  $y_k = q_k$ for all $k$ and 
 $- y_r \delta = q_{r, 1}$. Thus, $\delta = - \frac{q_{r,1}}{q_r}$, 
 and this number is positive. Indeed, the polynomial 
 $q$ consists of  $n+1$ exponents and has $n$ zeros, hence, 
all its coefficients  are nonzero and have alternating signs. Since  $q_{r,1}$ and $q_r$ 
are two consecutive coefficients, we have $\frac{q_{r,1}}{q_r} < 0$.  
Substituting to~(\ref{ineq1}) we obtain
 $$
 p_{\lambda}(t) \ =  \ p(t) \ + \ \lambda \, q(t)\ + \  O(\lambda^2)\, , \qquad \lambda \to 0\, .
 $$
 On the other hand, 
 $\|p +  q\| < 1$, and consequently  $\|p + \lambda q\| < 1 - c \lambda \, , \ \lambda \in (0,1]$. Therefore,   $\|p_{\lambda}\| < 1 - c \lambda$ for all positive  $\lambda$
 small enough. Moreover, since  
 $F(p+\lambda q) = F(p) + \lambda F(q)$ with $F(q) > 0$, it follows that  $F(p_{\lambda}) > F(p) + c_1\lambda$
for all positive  $\lambda$ small enough, where $c_1 > 0$ is a constant. 
This contradicts the optimality of~$p$.

   {\hfill $\Box$}
\bigskip

\begin{center}
\large{ \textbf{5. Corollaries and special cases}}
\end{center}
\bigskip

Theorem~\ref{th10} claims that the  maximal value of the $\ell$th 
derivative of a real quasipolynomial on the unit ball (in the monotone norm) 
strictly increases in each component~$h_i$. This, in particular, gives us the 
largest value of the $\ell$th derivative over all polynomials when the vector 
$\bh$ fills a rectangle in~$\re^n$.  
\begin{cor}\label{c10}
If the parameter $h_i$ belongs to the rectangle  $\{\bh \in \re^n_{++} \ | \  h_i \in (0, \alpha_i]\}$, then the maximal value of  $\Phi_{\ell} ({\bh})$ on the unit ball of an arbitrary monotone norm is achieved at a unique point, which is the vertex~$\bh = (\alpha_1, \ldots , \alpha_n)$.
\end{cor}

Let us now apply the comparison theorem to some concrete norms. 
We also call  a monotone norm on  $\re_+$   {\em shift-monotone}, 
if  $\bigl\|f(\cdot + a)\bigl|_{\re_+}\bigr\| \, \le \, \|f(\cdot)\|$ 
for every positive   $a$.
All the $L_p$-norms possess this property. 
For an arbitrary shift-monotone norm~$\|\cdot \|$, 
the maximal value of 
$\|p^{(\ell)}\|_{\infty}$ under the constraint~$\|p\|\le 1$ is equal to  $|p^{(\ell)}(0)|$, 
i.e., it is attained at zero. To prove this is suffices to note that this function is coercive 
and hence attains its maximum at some point and then 
shift  that point to zero. Thus,
\begin{cor}\label{c20}
Suppose that the monotone norm $\|\cdot \|$ is also a shift-monotone;  
then for every~$\ell \in \n$,  the maximal value of $\|p^{(\ell)}\|_{\infty}$ on the unit ball $p \in \RP,\ \|p\| \le 1$, is an increasing function in~$\bh\, $.
\end{cor}
As we know, the~$L_{\infty}$-norm of a polynomial over an arbitrary Chebyshev system
takes the largest value~$p^{(\ell)}(0)$ 
 on the ball $\|p\|_{\infty} \le 1$ at a unique polynomial, which possesses an  alternanse and  is called the 
 Chebyshev polynomial for that system. In our case this is the Chebyshev polynomial 
 on the  system~$E_{\bh}$ that has an alternance on~$[0, +\infty)$ of $n$ points including 
the point~$0$. We denote this polynomial by
$T_{\bh}$. For every $\bh \in \re^n_{++}$, the polynomial 
$T_{\bh}$ can be found numerically by the Remez algorithm~\cite{R1}, which makes it possible to compute~$\Phi_{\ell}(\bh)\, = \, m_{\ell}(\bh)$. 
We are  interested in the uniform sharp constant~$\, m_{\,\ell , n}$.  

\bigskip

\begin{center}
\large{ \textbf{6. The Markov-Bernstein inequality for exponential polynomials}}
\end{center}
\smallskip

\bigskip

For an arbitrary vector $\bh \in \co_{++}$ and for every~$\ell \ge 0$, the value
$\Phi_{\ell}(\bh)$ of the problem~(\ref{markov2}) for the functional  $F(p) = - p^{(\ell)}(0)$, is also the sharp constant in the Markov-Bernstein inequality for exponential polynomials
over  the system~$E_{\bh}$: 
\begin{equation}\label{mb}
\bigl\| p^{(\ell)} \bigr\|_{\infty} \quad \le \quad \Phi_{\ell}(\bh) \, \bigl\| p \bigr\|\, . 
\end{equation}
Let us note that each shift-monotone norm~$\|\cdot \|$
on $\re_{+}$ has its own constant $\Phi_{\ell}(\bh)$, 
therefore, the prober notation for it is~$\Phi_{\ell}(\bh, \, \|\cdot\|)$. 
We, however, 
use the previous short notation.  If~$\bh$ is a real vector
with components from the unit cube  $\cI_{n}\, = \, 
\{\bh \in \re^n, \ 0< h_i \le 1, \, i = 1, \ldots , n\}$, then the 
maximal value of the constant~$\Phi_{\ell}(\bh)$, according  to Corollary~\ref{c10}, 
is attained if~$\bh = \be = (1, \ldots, 1)$.  
Therefore, the maximum is always achieved at the polynomial $\cP_{\be}$, i.e., 
at a function of the form 
 $p(t)\, = \, e^{-t} q(t)$,  where $q(t)$ is an arbitrary 
 algebraic polynomial of degree~$n-1$.  In particular, 
for $\|\cdot \| = \| \cdot \|_{\infty}$, the extremal polynomial is~$p(t)\, = \, e^{-t} R_{n-1}(t)$, 
where~$R_{n-1}$ is a Chebyshev algebraic polynomial of degree~$n-1$ with the Laguerre
weight~$e^{-t}$, for which 
the function $e^{-t} R_{n-1}(t)$ possesses  $n$ points of alternance on 
the positive half-line~\cite{CLM, Sz}. 
Let us collect all these facts in the following theorem: 
\begin{theorem}\label{th25}
For an arbitrary shift-invariant norm 
$\|\cdot \|$ on $\re_+$, arbitrary~$\ell \ge 1$ and $\bh \in \cI_n$, we have 
\begin{equation}\label{mb1}
\bigl\| p^{(\ell)} \bigr\|_{\infty} \quad \le \quad \Phi_{\ell}(\be) \, \bigl\| p \bigr\|\ , 
\qquad p \in \RP_{\bh}\, .  
\end{equation}
This inequality becomes an equality at a unique, up to multiplication by a constant,
polynomial~$\RP_{\be}$. For $\|\cdot\| = \|\cdot \|_{\infty}$, one has $\Phi_{\ell}(\be) = m_{\ell, n}$, and a unique extremal polynomial 
is the exponential Chebyshev polynomial~$T_{\be} = e^{-t}R_{n-1}(t)$, 
where~$R_{n-1}$  is an algebraic Chebyshev polynomial of degree~$n-1$ 
with the Laguerre weight~$e^{-t}$.
Moreover,  $m_{\ell, n} \, = \, |T_{\be}^{(\ell)}(0)|$. 
\end{theorem}
Thus,~$m_{\ell, n}$ is equal to the 
the $\ell$-th derivative of the polynomial  $e^{-t}R_{n-1}(t)$ at zero. 
For every~$n$, this polynomial can be explicitly found, hence, the value~$m_{\ell, n}$ 
is efficiently computable. However, a natural question arises 
how fast does it grow in~$n$ and in $\ell$. 
This problem is reduced to the estimation of the  $\ell$-th derivative of the polynomial 
$R_{n-1}$ at zero. The first upper bound was presented in 1964 by 
Sege~\cite{Sz}, who proved that  $|R_{n-1}'(0)| \le Cn$. Then this result was successively 
improved in~\cite{CLM, Fre,  MN}. Asymptotically sharp estimates  
for all $\ell$ have been obtained by V.Sklyarov~\cite{S}:
\smallskip

 \noindent \textbf{Theorem A} (V.Sklyarov, 2009). {\em 
For arbitrary natural~$n\ge 2$, let  $R_{n-1}$ denote the algebraic 
Chebyshev polynomial of degree~$n-1$ with the Laguerre weight~$e^{-t}$. 
Then for each $\ell \in \n$, we have 
\begin{equation}\label{skl}
\frac{8^{\ell}(n-1)\, !\, {\ell}\, !}{(n-1-{\ell})\, !  \, (2{\ell})\, !} \Bigl(1 - \frac{{\ell}}{2(n-1)} \Bigr)\quad \le \quad \Bigl|\, R_{n-1}^{(\ell)}(0)\, \Bigr|\quad \le \quad \frac{8^{{\ell}}(n-1)\, !\, {\ell}\, !}{(n-1-{\ell})\, !  \, (2{\ell})\, !}.
\end{equation}
}
\medskip
\smallskip

Thus, the ratio between the upper and the lower bound~(\ref{skl}) \
tends to one as $n\to \infty$ when~$\ell$ is fixed. 
To  estimate~$m_{\, \ell, n}$ one needs to evaluate the derivatives of the polynomial~$T_{\be}$. 
They can easily be obtained with the derivatives of the polynomial~$R_{n-1}$: 
$$
\Bigl[\, e^{-t}R_{n-1}(t)\, \Bigr]^{(\ell)}\quad = \quad 
e^{-t} \, \sum_{j=0}^{\ell} \  (-1)^{\ell-j}\, {k \choose j}\ R_{n-1}^{(j)}(t)\, .
$$  
Taking into account that the signs of~$R_{n-1}^{(j)}(0)$ alternate in~$j$, we get  
\begin{equation}\label{LM}
m_{\, \ell, n} \quad = \quad 1 \ + \ \sum_{j=1}^{\ell} \, \, {\ell \choose j}\  \, 
\bigl|\, R_{n-1}^{(\ell)}(0)\, \bigr|.
\end{equation}
Now invoke~(\ref{skl}) and obtain after elementary simplifications: 
\begin{equation}\label{LM1}
\sum_{j=0}^{\ell} \, \,  \left(1 - \frac{{j}}{2(n-1)} \right)
\frac{8^j \ {n-1 \choose j} \ {\ell \choose j}}{{2j \choose j}}\quad \le \quad  
m_{\ell, n} \quad \le \quad  \sum_{j=0}^{\ell} \, \, \frac{8^j \ {n-1 \choose j} \ {\ell \choose j}}{{2j \choose j}} \,
\end{equation}
(all the terms with $j \ge d$ are zeros;  ${n \choose 0} = 1$ for each $n \ge 0$). 
We have not succeeded in simplifying these expressions. Combining with 
Theorem~\ref{th25}, we obtain
\begin{theorem}\label{th35}
For all $\ell, n\in \n, \, n \ge 2,$ and $\bh \in \cI_{n}$, we have 
\begin{equation}\label{markov25}
\max_{p \in \RP(\bh)}\ \bigl\|p^{(\ell)}\bigr\|_{\infty} \quad \le \quad  m_{\, \ell , n}\, \bigl\|p \bigr\|_{\infty}\, , 
\end{equation}
where the equality is attained at the polynomial~$T_{\be} = e^{-t}R_{n-1}$, and 
\begin{equation}\label{markov3}
\sum_{j=0}^{\ell} \, \,  \left(1 - \frac{{j}}{2(n-1)} \right)
\frac{8^j \ {n-1 \choose j} \ {\ell \choose j}}{{2j \choose j}}\quad \le \quad  
m_{\, \ell , n} \quad \le \quad  \sum_{j=0}^{\ell} \, \, \frac{8^j {n-1 \choose j} {\ell \choose j}}{{2j \choose j}}
\end{equation}
(all the terms with  $j\ge n$ are zeros).
\end{theorem}
\begin{cor}\label{c40}
Under the assumptions of Theorem~\ref{th25},
if the half-interval $(0,1]$ for  $h_i$ is replaced by $(0, \alpha]$,
then the estimate~(\ref{markov3}) is replaced by $\alpha^{\ell}$.
This inequality is asymptotically sharp with the extremal polynomials~$e^{-t/\alpha }
R_n\bigl(t/\alpha\bigr)$.
\end{cor}

Thus, for every $\ell, n$, we have an upper bound  for the sharp constant~$m_{\ell, n}$ 
in the Markov-Bernstein inequality for real exponential polynomials. 
In the next section we deal with applications to dynamical systems, where 
we need only the case $\ell=2$, for which the inequality~(\ref{markov3}) 
gets the following form: 
\begin{equation}\label{k2}
m_{2, n} \ \le  \ \frac{16 \, n^{\, 2} \, - \, 24 \, n \, + \, 11}{3}\, .
\end{equation}
This estimate is sharp 
asymptotically as~$n\to \infty$. For small~$n$, it can be computed precisely.  
Table~\ref{tab1} presents the values~$m_{2, n}$ for $n\le 10$. 
We see that the general  estimate~(\ref{th35})
in the third column is sufficiently close to the sharp estimate~(column~2). 
\bigskip 

\begin{table}[thb]
\begin{center}
\begin{tabular}{|c|c|c|}\hline
 $n$ & $\mbox{\rm the upper bound for}\ m_{\,2, n}$ 
 & $\mbox{\rm the bound~(\ref{k2})}$\\
 \hline
\rule{0pt}{9pt}\noindent
 $2$ &  8.182 & 9 $\ {}$\\
  \hline
 $3$ &  25.157 & 27$\frac23$\\
  \hline
 $4$ &  52.587 & 57 $\ {}$ \\
  \hline
 $5$ &  90.585 & 97  $\ {}$\\
  \hline
 $6$ &  139.191 & 147$\frac23$\\
  \hline
 $7$ & 198.420 & 209  $\ {}$ \\
  \hline
 $8$ &  268.283 & 281  $\ {}$ \\
  \hline
 $9$ &  348.788 & 363$\frac23$\\
  \hline
 $10$ & 439.938 & 457  $\ {}$ \\
 \hline
\end{tabular} \\[2mm]
\caption{\footnotesize{The upper bound on~$m_{2, n}$ for $n = 2, \ldots , 10.$}}\label{tab1}
\end{center}
\end{table}
\smallskip

\begin{remark}\label{r5}
{\em  Since $M_{\ell, n}\, \ge \, m_{\ell, n}$, 
the upper bound on~$m_{\ell, n}$ does not necessarily give an upper bound for~$M_{\ell, n}$. 
In general, Theorem~\ref{th35} is not valid for complex exponential polynomials. 
In particular, we do not know whether the estimates~(\ref{markov25}) 
give also estimates for~$M_{\ell, n}$. In  Conjecture~\ref{conj10}
we suppose that the answer is affirmative and   $M_{\ell, n}\, = \, m_{\ell, n}$. }
\end{remark}

\bigskip

\begin{center}
\large{ \textbf{7. The stability of linear switching systems}}
\end{center}
\smallskip

\bigskip 

Linear switching system is a linear ODE  $\dot \bx (t) \, = \,  A(t)\bx(t)$
on the vector-function $\bx \, : \, \re_+ \, \to \, \re^n$ with the initial 
condition~$\bx(0) = \bx_0$ with a matrix function $A(t)$ 
taking values from a given compact set~$\cA$ called 
 {\em control set}. The {\em control function}, 
 or the  
{\em switching law}  is an arbitrary measurable function  
$A: \re_+ \to \cA$. Linear switching systems naturally appear  
in problems of robotics, electronic engineering, mechanics, planning etc.~\cite{L}.
One of the main issues  is to find or estimate the fastest possible growth of trajectories
of the stability of the system. The {\em Lyapunov exponent} 
$\sigma (\cA)$ of the system is the infimum of the numbers  $\alpha$ for which  
    $\|\bx(t)\| \, \le \, C\, e^{\alpha  t} \|\bx_0\|, \ \ t \in [0, +\infty)$. 
We shall identify the linear switching system with the corresponding 
matrix family (control set)~$\cA$.  The Lyapunov exponent
does not change after replacing the control set by its convex hull. 
Therefore, without lost of  generality we assume that~$\cA$ is convex. 
Moreover, it can also be assumed that the system is irreducible, i.e., its matrices 
do not share common invariant nontrivial subspaces.

A system is called {\em asymptotically stable} (in short,  ``stable''), if all its trajectories tend to zero. 
If $\sigma < 0$, then the system is obviously stable. The converse is 
less obvious: the stability implies that~$\sigma < 0$~\cite{MP2}.
For one-element control sets~$\cA = \{A\}$, the stability is equivalent to that 
the matrix  $A$ is Hurwitz, i.e., all its 
eigenvalues have  negative real parts.  
If $\cA$ contains more than one matrix, the stability problem becomes 
much harder.  It is well known~\cite{MP1, O} that the stability 
is equivalent to the existence of 
    {\em Lyapunov norm} in~$\re^n$, i.e., a norm for which there exists~$\delta > 0$ 
    such that
     $\|\bx(t)\| \, \le \, e^{- \, \delta t}\,  \|\bx_0\|\, , \ t \in \re_+$ for all  
    trajectories~$\bx(t)$. N.Barabanov~\cite{B} showed that
   for an arbitrary convex control set~$\cA$, there exists an 
   {\em invariant Lyapunov norm} such that  
    $\, \|\bx(t)\| \, \le  \, e^{\, \sigma\, t}\,  \|\bx_0\|$ for every  
    trajectory~$\bx(t)$, and for every point~$\bx_0$ there exists a trajectory
     $\tilde \bx(t)$ such that  $\tilde \bx(0) \, = \, \bx_0$ and  $\|\tilde\bx(t) \| \, = \, e^{\, \sigma\, t}\,  \|\bx_0\|$.

The problem of computation of the Lyapunov exponent of the system is equivalent to the solution of the 
stability problem. Indeed, for an arbitrary~$s$, the inequality $\sigma(\cA) < s$
is equivalent to that $\sigma(\cA - sI) < 0$, i.e., is equivalent to the stability of the 
system  $\cA - sI $, where  $I$ is a unit matrix and  $\cA - sI  \, = \, \{A - sI \ | \ A\in \cA\}$. Consequently, if we are able to solve the stability problem for every matrix 
family, then we can compare the Lyapunov exponent with an arbitrary number, 
which allows us to compute the Lyapunov exponent merely by the double division.
Unfortunately, the  general stability problem is quite 
difficult; the existing methods either work in low dimensions 
(as a rule, at most~$4-5$), or give too rough estimates. 
For example the method of common quadratic Lyapunov function (CQLF) 
gives only necessary conditions for stability, 
which is far from being necessary~\cite{L, LM}. 
Other methods, for example, by piecewise-linear or by piecewise-quadratic 
Lyapunov functions~\cite{BM1, BM2}, the extremal polytope method~\cite{GLP17}, etc., 
are sharper but, as a rule, are realizable 
only in small dimensions. The discretization method considered below is well known and 
looks quite prospective. 

\bigskip 

\begin{center}
 \textbf{7.1.  Discrete vs continuous}
\end{center}
\smallskip

\bigskip

The disctetisation method reduces the stability of the linear switching system 
to the stability of a suitable discrete-time system: 
 $\bx(k+1) = B(k)\bx(k), \ B(k) \in\cB, \ k \ge 0$, where $\cB$ 
 is a compact matrix set, see~\cite{BCM, BS, PJ1, PJ2}.  
The following key statement has been established in~\cite{MP2}: 
\smallskip

\noindent \textbf{Fact 1.} {\em If for some~$\tau_0 > 0$, the discrete 
system generated by the family of matrices  $\cB = I + \tau_0 \cA$ is stable, 
then it is stable for all
$\tau < \tau_0$, and, moreover, so is the continuous time system~$\cA$.}

\smallskip 

Since 
$e^{\tau A} \, = \, I + \tau A + O(\tau^2)$,  the matrix  $B = I + \tau A$
 gives a linear approximation of $e^{\tau A}$. 
In fact, we approximate every trajectory of the system~$\cA$ 
by the Euler method with the step size~$\tau$ and analyse the 
corresponding piecewise-linear trajectories. If the step size 
tends to zero, then, obviously, the stability of piecewise-linear approximation 
implies the stability of the original system. 
Fact~1, however, asserts more: if the piecewise-linear approximation is stable 
for some step size~$\tau_0$, not necessarily small, 
then it is stable for every step size and the system~$\cA$ is stable. 

Thus, the Euler method with the step~$\tau$ not only approximates the trajectories of the 
system, but also gives sufficient conditions for stability.
The main issue is to what extend that condition is necessary?  The fundamental 
problem is formulated as follows: 

\smallskip

\noindent \textbf{Problem 1.} {\em For a given continuous-time system~$\cA$ and for arbitrary $\varepsilon > 0$, find $\tau_0 > 0$ such that the inequality $\sigma(\cA)< - \varepsilon$ 
implies instability of the discrete-time system  
$\cB = I + \tau \cA$ for all $\tau \le \tau_0$. 
}

\smallskip 

In other words, how small the discretisation step~$\tau$ should be  to guarantee that 
the stability of the obtained discrete system implies the stability of the
continuous-time system with precision~$\varepsilon$~? 
The stability of the discrete system~$\cB$ is decided in terms of  
{\em the joint spectral radius}:  
$$
\rho(\cB)\quad = \quad \lim_{k \to \infty}\ \max_{B(j) \in \cB, j = 1, \ldots , k}\ \bigl\|B(k)\cdots B(1)\bigr\|^{1/k}\, . 
$$
A discrete-time system~$\cB$ is stable precisely when~$\rho(\cB) < 1$, see.~\cite{MP1}. 
Thus, if for some  $\tau_0$, we have  $\rho(I + \tau_0 \cA) < 1$, then 
$\rho(I + \tau \cA) < 1$ for all $\tau < \tau_0$ and, moreover,~$\sigma(\cA) < 0$.  
The inequality~$\rho(\cB)< 1$ is equivalent to 
the existence of the norm in~$\re^n$ in which $\|B\| < 1, \, B\in \cB$, 
i.e., the operators~$B\subset \cB$ all contractions in the corresponding norm.  
Considering the unit ball in that norm gives an equivalent condition: 
there exists a symmetrized convex bodies $G\in \re^n$ for which  
$B\, (G)\,  \subset \, {\rm int}\, G$ for all~$B \in \cB$.  
For a one-matrix family~$\cB$, the joint spectral radius becomes 
the usual spectral radius, i.e., the largest modulus of eigenvalues of the matrix.  
If $\cB$ contains at least two matrices, then the computation of its joint spectral radius 
becomes an   NP-hard problem~\cite{BT}. Nevertheless, recently 
several efficient methods of computation of $\rho(\cB)$ 
for generic matrix families were presented. They work well in dimensions $n \le 25$, 
and for nonnegative families, the dimension can be increased to several thousands~\cite{GP13, Mej, MR14}. Hence,  the stability problem for discrete-time systems can be efficiently solved
by computing the joint spectral radius. The solution of the Problem~1 makes it possible
to extend this method to the continuous-time system.  
To this end, the discretisation step~$\tau$ should not be too small. 
 Otherwise, all the matrices of the family $I + \tau \cA$ 
 will be close to the identity matrix, in which case all known methods of the 
 joint spectral radius computation suffer.

\begin{defi}\label{d30}
For a given compact family of matrices,   $\cA$ 
and for arbitrary $\varepsilon > 0$, let 
\begin{equation}\label{Kind}
S_{\cA}(\varepsilon) \ = \
\left\{
\begin{array}{lcl}
+\infty & \ , \ & \sigma(\cA) \, \ge \, - \varepsilon\, , \\
\sup\, \bigl\{\tau > 0 \ | \ \rho(I + \tau \cA) < 1\, \bigr\} & \ , \ & \sigma(\cA) \, < \, - \varepsilon\, .
\end{array}
\right.
\end{equation}
\end{defi}
We use the simplified notation~$S_{\cA}(\varepsilon) = S$. 
Thus, if   $\sigma(\cA) < - \varepsilon$, then 
 $\rho(I + \tau \cA) < 1$ for all  $\tau \, < \, S$.
The value~$S$ has the following meaning:  if we do not know the value of $\sigma (\cA)$, 
but have an upper bound for  $S$, then we choose an arbitrary value 
 $\tau < S$ and compute the joint spectral radius of the operator $\rho(I + \tau \cA)$. 
If it is larger than or equal to  one, then $\sigma \ge - \varepsilon$;
otherwise   $\sigma < 0$. Denote also  
$$
S_r(\varepsilon)\ = \ \inf\, \Bigl\{ \ S_{\cA} (\varepsilon) \  \Bigl| \   \max_{A \in \cA}\rho(A)\, \le \, r  \ \Bigr\}\, .
$$
Thus, $S_r(\varepsilon)$ is a lower bound for  $S_{\cA}(\varepsilon)$ which is 
uniform over all families~$\cA$ of~$n$ matrices with the spectral radii not exceeding~$r$. 
In Theorems~\ref{th10} and \ref{th5} we express~$S_{\cA}(\varepsilon)$ and $S_r(\varepsilon)$ 
by the values~$M_{2}(\bh)$ and 
$M_{2, n}$ for $\bh \, = \, -  \, {\rm sp}\, (A), \ 
A\in \cA$.

Now our  plan is to estimate the value $S_{\cA}(\varepsilon)$ from below by applying the Lyapunov norm
of the family~$\cA$ (Propositions~\ref{p10} and~\ref{p30}). 
After this the problem will be reduced to the solution of extremal problem~(\ref{prob1}), 
which is equivalent to finding the sharp constant in the Markow-Bernstein type 
inequality for exponential polynomials
(Proposition~\ref{p20}). Then deriving a uniform bound for~$S_{r}(\varepsilon)$ 
will be equivalent to the minimization of the value~$S_{\cA}(\varepsilon)$ 
over all families of matrices~$\cA$
 with spectral radii at most~$r$. This will be done in Theorem~\ref{th5}.
\smallskip

\bigskip

\begin{center}
\textbf{7.2. Reduction to an extremal problem}
\end{center}
\smallskip

\bigskip 

We are going to estimate the discretization step~$S_{\cA}(\varepsilon)$ from below 
by solving  a special extremal problem. 
Let $\bh = (h_1, \ldots , h_n) \in \co^n$ and 
$\RP_{\,\bh}$ be the corresponding space of real quasipolynomials.
For  given~$\varepsilon > 0$,  we denote by~$s\, (\bh, \varepsilon)$ 
the value of the following problem:
 \begin{equation}\label{prob1}
 \left\{
 \begin{array}{l}
 \frac{1 - p(0)}{p'(0) - \varepsilon p(0)} \ \to \ \min , \\
 {}\\
 p \in {\RP}_{h}\\
  \|p\|_{\re_+} \le 1\, , \ p\, '(0) \, > \, \varepsilon p(0)\, .
 \end{array}
 \right.
 \end{equation}
The geometric sense of this problem becomes obvious by Proposition~\ref{p10}. 
To formulate it we need to introduce some further notation. For a Hutwitz 
matrix~$B$ and for a vector $\bx \in \re^n$, we set $G_B(\bx) =
{\rm co}\, \{\, \pm \, e^{\, t\, B}\bx\, , \ t \in [0, +\infty)\}$. 
Thus, $G_B(\bx)$ is the convex hull 
of the curve  
$\{\, e^{\, t\, B}\bx\, , \ t \in [0, +\infty)\, \bigr\}\, \bigcup \, 
\{\, - \, e^{\, t\, B}\bx\, , \ t \in [0, +\infty)\, \bigr\}$, which connects the points 
$\bx$ and $\, - \bx$. 

\begin{prop}\label{p10}
Let $B$ be a Hurwitz matrix,   $\bx \in \re^n \setminus \{0\}$ is an arbitrary point and   $\varepsilon > 0$; then the largest number $\tau$ for which   $\bx + \tau \, (B\, - \, \varepsilon I)\, \bx \, \in \, G_{B}(\bx)$ is equal to 
$s\, (\bh, \varepsilon)$, where $\, \bh \, = \, -{\rm sp}(B)\ $ (Fig. 1)
\end{prop}
\begin{center}
\begin{figure}[h!]
\centering
 	{\includegraphics[scale = 0.4]{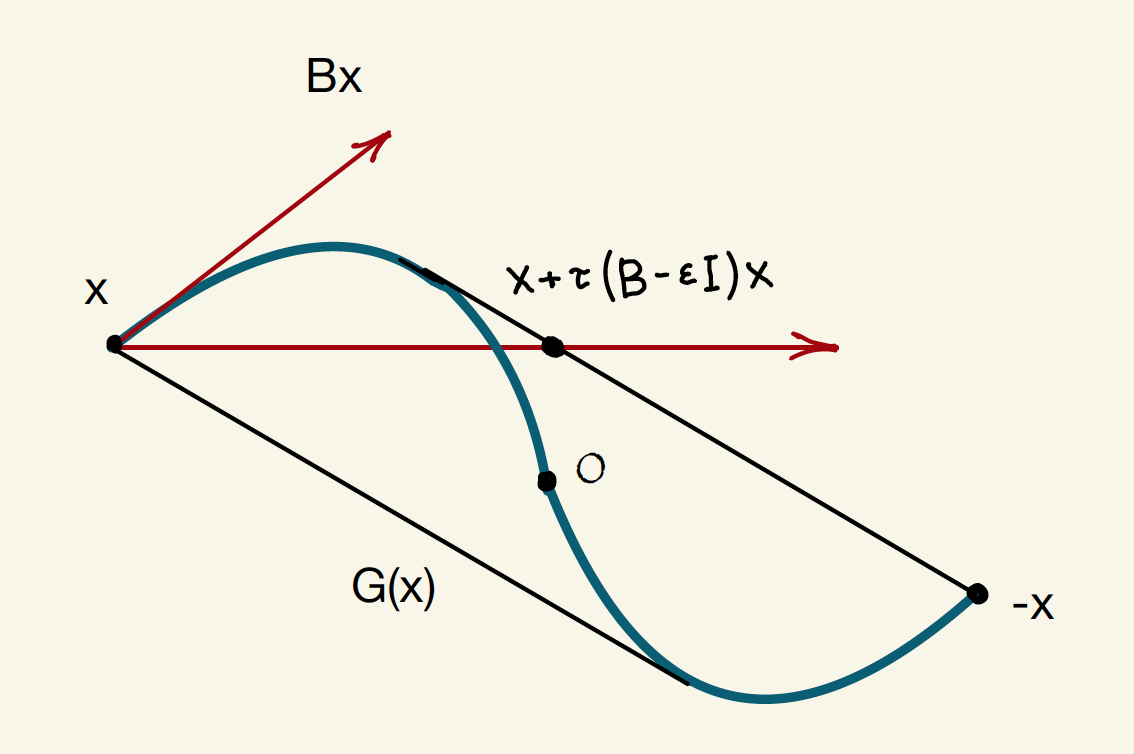}}
 	\caption{The largest number $\tau$ for which   $\bx + \tau \, (B\, - \, \varepsilon I)\, \bx \, \in \, G(\bx)$}
 	\label{fig. 1}
 \end{figure}
 \end{center}

{\tt Proof.} Denote  $G_B(\bx) = G$. Assume that  $B$ possesses 
 $ n$ different eigenvalues; the general case then follows by the limit passage. 
 By the Caratheodory theorem,  an arbitrary point of the set~$G$ 
 is a convex combination of at most  
     $n+1$ its extreme points of the form,  i.e., points of the form 
      $\pm \, e^{\, t\, B}\bx\, , \ t \ge 0$.
Hence, there are   $n+1$ nonnegative numbers  
 $\{t_k\}_{k=1}^{n+1}$ and numbers $\{q_k\}_{k=1}^{n+1}$
such that 
\begin{equation}\label{eq10}
\bx \ + \ \tau \, (B\, - \, \varepsilon I)\, \bx \quad = \quad \sum_{k=1}^{n+1} \, q_k \, e^{\, t_k\, B}\, \bx\, , \qquad
\sum_{k=1}^{n+1} |q_k| = 1\, .
\end{equation}
In the basis of eigenvectors of~$B$,  
we obtain the diagonal matrix 
$\tilde B \, = \, {\rm diag} \, (-h_1, \ldots , -h_n)\, $ and the point  
$\tilde \bx \, = \, \bigl(\tilde x_1, \ldots , \tilde x_n  \bigr) \in \co^n$. 
Moreover,  
 $\, {\rm Re\, } \, h_k > 0\, , k = 1, \ldots , n$.
We assume that all coordinates of the vector  $\tilde x_k$ are nonzero; 
the general case will then follow by the limit passage. 
Equality~(\ref{eq10}) has the following form  in the new basis:  
$$
\Bigl(\, I \ - \ \tau \,  {\rm diag} \, \bigl(h_1 + \varepsilon, \ldots , h_n + 
\varepsilon)\ \Bigr)\, \tilde \bx \,   \ = \  \sum_{k=1}^{n+1} \, q_k \, e^{\, - t_k\, h_j}\tilde \bx\, .
$$
For each component  $\tilde x_j$, we obtain: 
$$
\bigl(\, 1 \, - \, \tau \, (h_j \, + \, \varepsilon\, )\, \bigr)\, \tilde x_j \,   \ = \  \sum_{k=1}^{n+1} \, q_k \, e^{\, - t_k\, h_j}\tilde x_j\, , 
$$
then we divide by  $\tilde x_j$:
\begin{equation}\label{eq14}
 1 \, - \, \tau \, (h_j \, + \, \varepsilon\, )  \quad = \quad  \sum_{k=1}^{n+1} \, q_k \, e^{\, - t_k\, h_j}\, \quad j = 1, \ldots , n\, . 
\end{equation}
This equation does not contain coordinates of the vector~$\tilde \bx$.   
The maximal~$\tau$ for which there are coefficients  $\{t_k\}_{k=1}^{n+1}$ and 
$\{q_k\}_{k=1}^{n+1}$ such that~(\ref{eq14}) holds, must be the same for all 
 $\tilde \bx \ne 0$. We set  $\tilde \bx = \be$
(the vector of ones), and denote by $\tilde G$ 
the set of all values of the sum~$\sum_{k=1}^{n+1} \, q_k \, e^{\, - t_k\tilde B}\be$
for all suitable~$\{t_k\}_{k=1}^{n+1}$ and $\{q_k\}_{k=1}^{n+1}$. 
This is a convex body in~$\co^n$ over the field of real numbers, i.e., 
 a convex body in~$\re^{2n}$.  By the convex separation theorem, for an arbitrary $\delta > 0$,  the point  $\be \ + \ \tau \, (\tilde B\, - \, \delta I)\, \be$ does not belong to~$\tilde G$  precisely when it can be strictly separated from it with a linear functional, i.e., when there exists a nonzero vector 
 $\bz = (z_1, \ldots , z_n) \in \co^n$, for which 
\begin{equation}\label{eq16}
 {\rm Re\, } \ \Bigl( \, \be \, + \, \delta\, (\tilde B \, - \, \varepsilon I 
 \, )\, \be\, , \, \bz \Bigr) \quad > \quad \max_{t \ge 0} \ \left| \, 
 {\rm Re\, } \left( \ \sum_{k=1}^n  \, \bar z_j \, e^{\, - t \, h_j } \, \right)\, \right| \, .
\end{equation}
Denote by   $\, q(t) \, = \, \sum_{j=1}^n\, \bar z_j e^{\, - \, t \, h_j } \in \cP_{\bh}$ a complex polynomial and by $p(t) \, = \, {\rm Re}\, q(t)$ the corresponding real quasipolynomial. Then the right-hand side of the inequality~(\ref{eq16}) is equal to~$\bigl\|\, p\, \bigr\|_{\re_+}$. 
Then we modify the left-hand side as follows: 
$$
\Bigl( \, \be  \, +  \, \delta (\tilde B \, - \, \varepsilon I)\, \be\, , \, \bz\, \Bigr)\  =\
\sum_{j=1}^n \bar z_j \, + \, \delta \, \sum_{j=1}^n (-h_j - \varepsilon)\, \bar z_j\  =
\  q(0) \, + \, \delta \,  q\, '(0) \, - \, \delta \, \varepsilon \, q(0)\, .
 $$
Thus, the left-hand side of~(\ref{eq16}) is equal to
$ p(0) + \delta \,  p'(0)  - \delta \, \varepsilon \, p(0)\, ,$
and we arrive at the inequality 
$$
 p(0) \, + \, \delta \,  p\, '(0) \, - \, \delta \, \varepsilon \, p(0)\quad 
 > \quad \|p\|_{\re_+}\, , 
$$
hence, 
$$
\delta \,  \bigl( \, p\, '(0) \, - \, \varepsilon \, p(0)\, \bigr) \quad 
 > \quad \|p\|_{\re_+}\, - \, p(0)\, . 
$$
The right-hand side is obviously nonnegative, consequently $p'(0) -  \varepsilon \, p(0) > 0$. 
Normalising the polynomial~$p$ so that $\|p\|_{\re_+} = 1$, we obtain 
$\delta  >  \frac{1 - p(0)}{p\, '(0) - \varepsilon p(0)}$.Thus, the point ${\be \, + \, \delta \, (B \, - \, \varepsilon I)\, \be}$ does not belong to the convex body
 $G_{B}(\be)$ if and only if there exists a convex body $p \in \RP_{\bh}, \, \|p\|_{\re_+} \le 1$ satisfying that inequality. 

{\hfill $\Box$}
\medskip

\begin{remark}\label{r10}
{\em In fact we prove more: if~$\tau < s(\bh, \varepsilon)$ and $\bx \ne 0$, then the point 
 $\bx + \tau \, (B\, - \, \varepsilon I)\, \bx$ lies inside~$G_B(\bx)$.}
\end{remark}

The lower bound for the estimate~$S_{\cA}(\varepsilon)$ 
is provided by the maximal possible value of~$s(\bh,  \varepsilon)$ 
over all vectors
  $\bh \, = \, -\, {\rm sp} \, (A) \, - \, \varepsilon \, \be, \ A\in \cA$.
We prove this in the following proposition by applying the invariant norm of the matrix family~$\cA$. 
 \begin{prop}\label{p30}
For an arbitrary matrix family $\cA$, we have 
\begin{equation}\label{sr}
\, S_{\cA} (\varepsilon) \quad \ge \quad  
\min_{A \in \cA}\, s\, (\bh, \varepsilon)\, , \
\end{equation}
$\mbox{where}\  \bh \ = \ -{\rm sp}(A) \, - \, \varepsilon\, \be\, .$
 \end{prop}
{\tt Proof.} It is needed to show that if $\sigma (\cA) < -\varepsilon$, then  
$\rho
 (I + \tau A) < 1$, for all~$\tau$ which are smaller than the right-hand side  of the inequality~(\ref{sr}).
Since   $\sigma (\cA) < -\varepsilon$, 
It follows that there exists a Lyapunov norm in~$\re^n$, 
for which   $\|\bx(t)\| \, < \, e^{-\varepsilon t} \|\bx_0\|, \ t>0, $ for every trajectory~$\bx(\cdot)$.
This is true, in particular, for trajectories without switches, i.e.,
for the stationary control. 
Thus, for an arbitrary matrix~$A \in \cA$, we have 
$\|e^{\, t A}\bx_0\| \, < \, e^{-\varepsilon t} \|\bx_0\|$, and therefore 
$\|e^{\, t\, (A + \varepsilon I)}\bx_0\|\, <  \, \|\bx_0\|$ for all $t > 0$. 
Hence, for each point   $\by$ from the symmetrised  convex hull of the set
 $\, \bigl\{\, e^{\, t\, (A + \varepsilon I)}\bx_0 \ \bigl| \ t \in \re_+\, \bigr\}$,
we have  $\|\by\| \le \|\bx_0\|$. Applying Proposition~\ref{p10} to the matrix~$B = A + \varepsilon I$ and taking into account Remark~\ref{r10},
we conclude that the point 
 $\, (I \, + \, \tau A) \bx_0$ belongs to the interior of the set $G_{A + \varepsilon I}(\bx_0)$,
and therefore, $\, \|(I \, + \, \tau A) \bx_0\| \, < \,  \|\bx_0\|$. Thus,
for each~$A \in \cA$, the operator norm of $I \, + \, \tau A$
is strictly smaller than one. Consequently,  $\rho (I + \tau \cA) <  1$.

{\hfill $\Box$}
\medskip

\begin{remark}\label{r21}
{\em The value~$s(\bh, \varepsilon)$
of the problem~(\ref{prob1}) is inversely proportional to its parameters: for every 
$\lambda > 0$, we have 
\begin{equation}\label{homog}
s \bigl(\lambda \bh\, , \, \lambda \varepsilon  \bigr) \ = \ \lambda^{-1}\, s \bigl(\bh, \varepsilon \bigr)\, .
\end{equation}
For the proof, it suffices to change the variable  $u = \lambda t$ and 
observe that 
 ${p_t}\,'(0) = \lambda \, {p_u}'(0)$.
Therefore, the computation or estimation of the value~$s (\bh, \varepsilon)$ is reduced to the case
$\bh \in \Delta$, i.e., $|h_n| \le 1$. 
}
\end{remark}

\smallskip

It remains to compute~$s(\bh, \varepsilon)$, i.e., to solve the problem~(\ref{prob1}). This can be done numerically  applying the convex optimization tools, since the set of real quasipolynomials 
~$q$ satisfying all the constraints of the problem is convex and the 
objective functions quasiconvex (all its level sets are convex). 
Therefore, for a concrete vector 
$\bh \in \co^n$, the problem can be efficiently solved, for instance, by the gradient relaxation methods. We, however, need a general lower bound for the value~$s\, (\bh, \varepsilon)$. 

\begin{prop}\label{p20}
For each   $\varepsilon > 0$ and $\, \bh \in \Delta$, we have 
 $$
s (\bh, \varepsilon) \quad > \quad \frac{2\, \varepsilon}{M_2(\bh)\, +\, 2\, \varepsilon^2}
 $$
where $s (\bh, \varepsilon)$ is the value of the problem~(\ref{prob1}).  
\end{prop}
{\tt Proof.} 
Let  us first forget for a moment that  $p$  is a quasipolynomial,
we leave the only assumption that~$\|p''\| \le \ M_2(\bh)\|p\|$. 
For the sake of simplicity we also denote~$M_2(\bh) \, = \, m$. 
Thus, we need to compute the minimum of the value~$\frac{1-p(0)}{p'(0) - \varepsilon p(0)}$ 
under the assumptions~$p \in C^2(\re_+), \, \|p\| \le 1$ and $\|p''\| < m$. 
Set $x = 1 - p(0)$, then the objective function in the problem~(\ref{prob1}) 
gets the form  $f(x) = \frac{x}{p'(0) - \varepsilon (1-x)}$. 
The numerator and the denominator are both positive, hence 
the minimum of this fraction is achieved  when~$p'(0)$ is maximal. 
Denote~$p'(0) = \alpha$. 
Since $p''(t) \ge - m$, we have for each~$t\in \re_{+}$: 
  $$
  p(t) \quad  \ge \quad  p(0) \, + \, p'(0)t \, - \, \frac{mt^2}{2}\quad = \quad 
  1\, -\, x \, + \, \alpha t \, - \, \frac{mt^2}{2}\, .  
  $$
On the other hand, $p(t) \le 1$ because $\|p\|_{\re_+}\le 1$. 
Therefore, the maximum of the quadratic function  
  $1\, -\, x \, + \, \alpha t \, - \, \frac{mt^2}{2}$  over $t \in \re_+$
 does not exceed~$1$. Computing the value at the vertex of the parabola, we obtain~$\frac{\alpha^2}{2m}\, + \, 1 - x \le 1$, and so 
$\alpha \le \sqrt{2mx}$. Hence, 
$$
f(x) \ \ge \ \frac{x}{\sqrt{2mx} - \varepsilon (1-x) }, . 
$$
The minimum of the left-hand side of this inequality under the positivity constraints 
for the numerator and denominator is attained at the point~$x_{\min} = \frac{2\varepsilon^2}{m}$, 
therefore 
$$
f(x) \ \ge \ f(x_{\min}) \ = \ \frac{2\varepsilon}{m+ 2\varepsilon^2}\, . 
$$
Thus, the minimum of the objective function  $\frac{1-p(0)}{p'(0) - \varepsilon p(0)}$
is equal to~$\frac{2\varepsilon}{m+ 2\varepsilon^2}$ and is achieved at the polynomial~$p$ defined above. Since the quadratic function does not belong to~$\RP_{\bh}$ and 
the unit ball in the space~$\RP_{\bh}$ is compact, it follows that the minimum of 
the objective function is larger than~$\frac{2\varepsilon}{m+ 2\varepsilon^2}$.

{\hfill $\Box$}
\medskip

Combining Propositions~\ref{p30} and \ref{p20} we obtain the estimate for the 
length o the discretization interval~$S_{\cA}(\varepsilon)$. 
\begin{theorem}\label{th3}
For every~$\varepsilon > 0$ and  $r > \varepsilon$, 
for every system of operators~$\cA$ such that  $\max_{A \in \cA} \rho(A)\, \le \, r$, 
we have 
\begin{equation}\label{individ}
S_{\cA}(\varepsilon) \quad \ge \quad \max_{A\in \cA} \ \frac{2\, \varepsilon}{M_2(\bh)\, +\, 2\, \varepsilon^2}\ , 
\end{equation}
where $\bh = \, -\, {\rm sp}\, (A) \, - \, \varepsilon \, \be$. 
\end{theorem}
The estimate~(\ref{individ}) is computed individually for all systems of matrices~$\cA$. 
A uniform estimate~$S_{r}(\varepsilon)$ over all families of matrices is given by the following  
\begin{theorem}\label{th5}
For every~$\varepsilon > 0$ and  $r > \varepsilon$, 
 for every system of operators~$\cA$ such that  $\max_{A \in \cA} \rho(A)\, \le \, r$, 
we have 
$$
S_{r}(\varepsilon) \ \ge \ \frac{2\, \varepsilon}{r^2M_{2, n}}
$$
\end{theorem}
{\tt Proof of Theorem~\ref{th5}.} 
Applying Propositions~\ref{p10} and~\ref{p30} we conclude that the value
$S_{\cA} (\varepsilon)$ is not smaller than the minimal value~$s(\bh, \varepsilon)$
over all vectors~$\bh = -{\rm sp}(A) - \varepsilon \be, \ A \in \cA$.  
Proposition~\ref{p20} yields  
$$
S_{\cA} (\varepsilon) \quad \ge \quad 
\min_{\bh = -{\rm sp}(A) - \varepsilon \be, \ A \in \cA}\quad 
\frac{2\, \varepsilon}{M_2(\bh)\, +\, 2\, \varepsilon^2} \, .  
$$
On the other hand 
$$
s(\bh, \varepsilon) \ \ge \ \frac{2\, \varepsilon}{(\rho(A) - \varepsilon)M_{2, n} \, +\, \frac{2\varepsilon^2}{\rho(A) - \varepsilon}}\, .
$$
If $\rho(A) \le r$ for all  $A \in \cA$, then  
$$
S_{r}(\varepsilon) \ \ge \ \frac{2\, \varepsilon}{(r - \varepsilon)^2M_{2, n}  \, +\, 2\, \varepsilon^2}
$$
Since   $M_{2, n}  > 2$ and $r \ge \varepsilon$, 
we see that the denominator of that fraction is at least 
$r^{\, 2} M_{2, n} $, which completes the proof.

{\hfill $\Box$}
\medskip

Theorems~\ref{th3} and \ref{th5} show that 
the lower bound for the discretization step is linear in~$\varepsilon$, 
which is much better than one might expect. As a rule, such estimates are exponential in~$\varepsilon$ and can hardly be applied in dimensions bigger that three. 
The estimate provided by Theorem~\ref{th5} depends on the dimension 
only in the coefficient of the linear function.  This makes it applicable even in relatively high dimensions. However, there is one difficulty on this way: we cannot compute~$M_{2, n}$. Only the families  $\cA$ that consist of matrices with real eigenvalues the constant~$M_{2, n}$ in Theorem~\ref{th5} becomes efficiently computable because it can be replaced by~$m_{2, n}$  
   (Theorem~\ref{th25}). 
\begin{cor}\label{c25}
If all the matrices of the family~$\cA$ have real spectra, then the value~$M_{2, n}$ 
in Theorem~\ref{th5} is replaced by~$m_{2, n}$ and  
the value~$M_{2}(\bh)$ in Theorem~\ref{th3} is replaced by~$m_{2}(\bh)$. 
\end{cor}
The value~$m_{2}(\bh)$ can be found with arbitrary precision by applying the  Remez algorithm. 
It is attained at the corresponding Chebyshev polynomial~$T_{\bh}$. 
The value~$m_{2, n}$, as it follows from Theorem~\ref{th25}, is attained at the polynomial~$e^{-t}R_{n-1}(t)$ and can be efficiently computed as well. It can be estimated by 
inequality~(\ref{markov3}) and, for small dimensions, can be evaluated precisely by 
computing  the Chebyshev polynomial~$R_{n-1}$
with the Laguerre weight. For small dimensions, the results are presented in Table~\ref{tab1}. We see that  $m_{2, n}$ grows in~$n$ not very fast, therefore, the estimate 
from Theorem~\ref{th5}, is applicable even for relatively high dimensions
 (few dozens). 

Certainly, the real spectrum   case is rather exceptional. 
In the general case one needs to compute or at east to estimate the value of  $M_{2, n}$
from above. We are not aware of any method of its computation with arbitrary precision. 
All known estimates are rough and can hardly be applicable. If Conjecture~\ref{conj10} is true, then   $M_{2, n} = m_{2, n}$,  and the problem will be solved. 
Otherwise, we have an open problem  of finding satisfactory  
upper bounds for~$M_{2, n}$. 

\vspace{1cm} 

The author expresses his sincere thanks to the Anonymous Referee for 
the attentive reading and many valuable remarks.

\bigskip


\bigskip

\begin{thebibliography}{}

\bibitem{B}
N.E.\,Barabanov,
\newblock {\em Absolute characteristic exponent of a class of linear nonstationary systems of differential equations},
\newblock Siberian Math. J. 29 (1988), 521--530.
\smallskip

\bibitem{BCM}
F.~Blanchini, D,~Casagrande and S.~Miani,
\newblock {\em Modal and transition dwell time computation in switching systems: a set-theoretic approach},
\newblock Automatica J. IFAC, 46 (2010), no 10, 1477--1482.

\bibitem{BM1}
F.\,Blanchini, S.\,Miani,
\newblock {\em Piecewise-linear functions in robust control},
\newblock Robust control via variable structure and Lyapunov techniques (Benevento, 1994),  213--243, Lecture Notes in Control and Inform. Sci., 217, Springer, London, 1996.
\smallskip



\bibitem{BM2}
F.\,Blanchini, S.\,Miani,
\newblock {\em A new class of universal Lyapunov functions for the control of
uncertain linear systems},
\newblock IEEE Trans. Automat. Control,   44  (1999), no 3, 641--647.
\smallskip

\bibitem{BT}
V.~Blondel and J.~Tsitsiklis,
\newblock {\em Approximating the spectral radius of sets of matrices in the max-algebra is
NP-hard},
\newblock IEEE Trans. Autom. Control, 45 (2000), No 9,  1762--1765.
\smallskip


\bibitem{BE0}
P.B.\,Borwein and T. Erd\'elyi
\newblock {\em Polynomials and polynomial inequalities},
\newblock Springer-Verlag, New. York, N.Y., 1995
\smallskip


\bibitem{BE1}
P.B.\,Borwein and T. Erd\'elyi
\newblock {\em Upper bounds for the derivative of exponential sums}
\newblock Proc. Amer. Math. Soc. 123 (1995), 1481 -- 1486.
\smallskip


\bibitem{BE2}
P.B.\,Borwein and T. Erd\'elyi
\newblock {\em A sharp Bernstein-type inequality for exponential sums},
\newblock J.Reine Angew. Math. 476 (1996), 127--141.
\smallskip


\bibitem{BE3}
P.B.\,Borwein and T. Erd\'elyi
\newblock {\em Newman's inequality for M\"untz polynomials on positive intervals},
\newblock J. Approx. Theory, 85 (1996), 132--139.
\smallskip

\bibitem{BS}
C.~Briat and A.~Seuret.
\newblock Affine characterizations of minimal and mode-dependent dwell-times
  for uncertain linear switched systems.
\newblock {\em IEEE Trans. Automat. Control}, 58(5):1304--1310, 2013.

\bibitem{CLM}
 H.\,Carley, X.\,Li, R.N.\,Mohapatra,
\newblock {\em A sharp inequality of Markov type for polynomials
associated with Laguerre weight},
\newblock J. Approx. Theory, 113 (2001), no 2, 221--228.
\smallskip


\bibitem{Dz}
V.K.\,Dzyadyk and I.A.\,Shevchuk,
\newblock  Theory of uniform approximation of functions by polynomials,
\newblock Walter de Gruyter, 2008.
\smallskip


\bibitem{Fr}
W.\,Fraser,
 \newblock {\em A Survey of methods of computing minimax and near-minimax polynomial approximations for functions of a single independent variable},
 \newblock J. ACM 12 (1965), no 295.
\smallskip

\bibitem{Fre}
G.\, Freud,
\newblock {\em  On two polynomial inequalities},
\newblock  Acta Math. Acad. Sci. Hungar., 22 (1971), no 1--2,
109--116.
\smallskip

\bibitem{G}
G.~Gripenberg.
\newblock Computing the joint spectral radius.
\newblock {\em Linear Algebra Appl.}, 234 (1996), 43--60.

\bibitem{GLP17}
N.~Guglielmi, L.~Laglia, and V.~Protasov.
\newblock Polytope {L}yapunov functions for stable and for stabilizable {LSS}.
\newblock {\em Found. Comput. Math.}, 17 (2017), no 2, 567--623.

\bibitem{GP13}
N.~Guglielmi and V.~Protasov.
\newblock Exact computation of joint spectral characteristics of linear
  operators.
\newblock {\em Found. Comput. Math.}, 13 (2013), no 1, 37--97.

\bibitem{K}
S.\,Karlin,
\newblock {\em  Representation theorems for positive functions},
\newblock  J. Math. Mech. 12 (1963), no 4,  599--618.
\smallskip

\bibitem{KS}
 S.\,Karlin and W.J.\,Studden,
 \newblock {\em Tchebycheff Systems: With Applications in Analysis and Statistics},
 \newblock Interscience, New York, 1996. 
\smallskip

\bibitem{KN}
M.G.\,Krein and A.A.\,Nudelman,
\newblock {\em The Markov moment problem and extremal problems: ideas and problems of P.L.Cebyshev and A.A.Markov and their further development},
\newblock Translations of mathematical monographs, Providence, R.I.,  v. 50 (1977).
\smallskip

\bibitem{L}
D.~Liberzon.
\newblock {\em Switching in systems and control}.
\newblock Systems \& Control: Foundations \& Applications. Birkh\"{a}user
  Boston, Inc., Boston, MA, 2003.

\bibitem{LM}
D.~Liberzon and A.~S. Morse.
\newblock Basic problems in stability and design of switched systems.
\newblock {\em IEEE Control Systems Magazine}, 19 (1999), 59--70.
\smallskip

\bibitem{Mej}
T.~Mejstrik,
\newblock {\em Improved invariant polytope algorithm and applications}, 
\newblock  ACM Trans.\ Math.\ Softw., 46 (2020), no  3 (29), 1--26.
\smallskip

\bibitem{MR14}
C.~M\"oller and U.~Reif,
\newblock {\em A tree-based approach to joint spectral radius determination},
\newblock Linear Alg. Appl., 563 (2014), 154-170.
\smallskip

\bibitem{MN}
L.\,Milev and N.\,Naidenov,
\newblock {\em Exact Markov inequalities for the Hermite and Laguerre
weights},
\newblock J. Approx. Theory, 138 (2006), no 1, 87--96. 
\smallskip

\bibitem{MP1}
A.P.\,Molchanov and E.S.\,Pyatnitskii,
\newblock {\em Lyapunov functions, defining necessary and sufficient
conditions for the absolute stability of nonlinear nonstationary control systems,}
\newblock Autom. Remote Control, 47 (1986), I -- no 3, 344--354, II - no 4, 443--451, III -- no 5, 620--630.
\smallskip

\bibitem{MP2}
A.~P. Molchanov and Y.~S. Pyatnitskiy.
\newblock Criteria of asymptotic stability of differential and difference
  inclusions encountered in control theory.
\newblock {\em Systems Control Lett.}, 13 (1989), no 1, 59--64.

\bibitem{New}
D.J.\,Newman,
\newblock {\em Derivative bounds for M\"untz polynomials},
\newblock J. Approx. Theory 18 (1976), 360--362.
\smallskip

\bibitem{O}
V.I.\,Opoitsev,
\newblock {\em Equilibrium and stability in models of collective behaviour},
\newblock Nauka, Moscow (1977).
\smallskip

\bibitem{PJ1}
V.Yu.Protasov and R.Jungers, 
 \newblock {\em  Is switching systems stability harder for continuous time systems?},
 \newblock  Proc. of 2013 IEEE 52nd Annual Conference on Decision and Control (CDC2013), Firenza (Italy), December 10--13, 2013.
\smallskip

\bibitem{PJ2}
V.Yu.Protasov and R.Jungers, 
\newblock {\em Analysing the stability of linear systems via
exponential Chebyshev polynomials},
\newblock  IEEE Trans. Automatic Control, 61 (2016), no 3,  795--798.
\smallskip

\bibitem{PJB}
V.\,Yu.~Protasov, R.\,M.~Jungers, and V.\,D.~Blondel,
\newblock {\em Joint spectral characteristics of matrices:
a conic programming approach},
\newblock SIAM J. Matrix Anal. Appl., 31 (2010), no 4,  2146--2162.
\smallskip

\bibitem{R1}
 E.Ya.\,Remez,
 \newblock {\em Sur le calcul effectiv des polynomes d'approximation des Tschebyscheff},
 \newblock Compt. Rend. Acade. Sc. 199, 337 (1934).
\smallskip

\bibitem{RS}
G.\,C.~Rota and G.~Strang,
\newblock {\em A note on the joint spectral radius},
\newblock Kon. Nederl. Acad. Wet. Proc. Vol. 63 (1960),  379--381.

\bibitem{S}
 V.P.~Sklyarov,
 \newblock {\em  The sharp constant in Markov's inequality for the Laguerre weight}, 
\newblock Sb. Math., 200 (2009), no 6, 887--897. 
\smallskip

\bibitem{Sz}
 G.\,Szeg\"o,
 \newblock {\em On some problems of approximations},
\newblock  Magyar Tud. Akad. Mat. Kutat\'o
Int. K\"ozl., 9 (1964), 3--9.
\smallskip


 
 

\end{thebibliography}
\end{document}